\documentclass[a4paper,11pt]{article}
\usepackage{amsfonts}
\usepackage{amssymb}
\usepackage{amsmath}
\usepackage{amsthm}
\usepackage{multirow}
\usepackage{diagbox}
\usepackage{amsopn}
\usepackage{gensymb} 
\usepackage{fancyhdr}
\usepackage{graphicx}
\usepackage{mathrsfs}
\usepackage{latexsym}
\usepackage{graphicx}
\usepackage{subfigure}
\usepackage{color}
 \usepackage{calc}
  \usepackage{tikz}
  
\usepackage{setspace}
\usepackage{textcomp} 

\usepackage[width=11cm]{geometry}   
\usepackage{varwidth}
\usepackage{algorithm}
\usepackage{algpseudocode}
\numberwithin{equation}{section} \allowdisplaybreaks
\theoremstyle{plain}

 \usepackage{rotating}

\theoremstyle{remark}

\numberwithin{equation}{section}

\renewcommand{\labelitemi}{$\bullet$}

\renewcommand{\labelitemiii}{$\diamond$}

\begin{document}
\title{\bf{An  ensemble Kalman filter approach based on level set parameterization for acoustic source identification using multiple frequency information}}

\author{\hspace{-2.5cm}{\small Zhi-Liang Deng$^{1}$ and Xiao-Mei Yang$^2$\thanks{
Corresponding author: yangxiaomath@163.com; yangxiaomath@home.swjtu.edu.cn
Supported by NSFC No. 11601067, 11771068 and No.11501087, the Fundamental Research Funds for the Central Universities No. 2682018ZT25 and ZYGX2018J085.
},} \\
\hspace{-1.5cm}{\scriptsize $1.$ School of Mathematical Sciences,  University of Electronic Science and Technology of China,
Chengdu 610054, China}\\
\hspace{-2.5cm}{\scriptsize $2.$ School of Mathematics,
Southwest Jiaotong University,
Chengdu 610031, China}
}
\date{}
\maketitle

\begin{abstract}
\noindent The spatial dependent unknown acoustic source is reconstructed according noisy multiple frequency data on a remote closed surface. Assume that the unknown function is supported on a bounded domain. To determine the support, we present a statistical inversion algorithm, which combines the ensemble Kalman filter approach with level set technique. Several numerical examples show that the proposed method give good numerical reconstruction.



\end{abstract}

\section{Introduction}
\begin{figure}\centering
 \begin{tikzpicture}
\draw [blue] (-2,-2) rectangle (4,4);
\draw [fill=yellow,ultra thick] (1,1) circle [radius=2];
\draw[fill=blue,thick] (0,0.4) to [out=90,in=145] (2,1.5);
\draw[fill=blue,thick] (2,1.5) to [out=-35,in=-95] (0,0.4);
\node at (1,1) {$D$};
\node at (3, 3) {$\Omega$};
\node at (2, 2) {$D_1$};
\draw 
(-2,-2) circle (2pt) 
(-1,-2) circle (2pt) 
(0,-2) circle (2pt) 
(1,-2) circle (2pt) 
(2,-2) circle (2pt) 
(3,-2) circle (2pt) 
(4,-2) circle (2pt)    
(4,-1) circle (2pt)   (4,0) circle (2pt) (4,1) circle (2pt) (4,2) circle (2pt) 
(4,3) circle (2pt) (4,4) circle (2pt) (3,4) circle (2pt) (2,4) circle (2pt)
(1,4) circle (2pt)  (0,4) circle (2pt) (-1,4) circle (2pt) (-2,4) circle (2pt) 
(-2,3) circle (2pt) (-2,2) circle (2pt) (-2,1) circle (2pt) (-2,0) circle (2pt) 
(-2,-1) circle (2pt)
 node[align=right,  below] {};
\end{tikzpicture}
\caption{The problem geometry.}\label{fig:digit}
\end{figure}
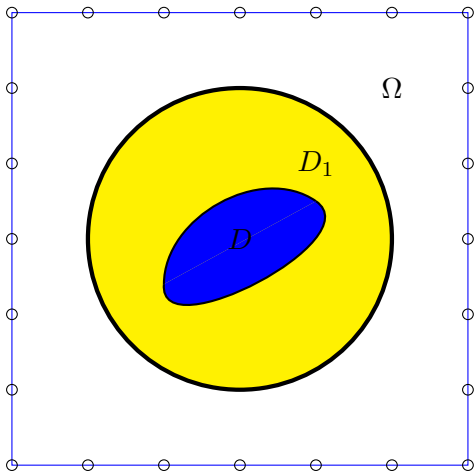

The unknown source problems occur widely in many real applications, e.g., remote sensing, radar detection \cite{kruk}, brain image etc. 
In this paper, we consider a reconstruction problem of acoustic sources from remote measurements of the acoustic field.  For an acoustic source with density $\mathcal{S}(k, x)$ supported on region $D\subset\subset\Omega$,  the pressure $u$ of the radiated  time-harmonic wave can be modeled \cite{eller} by 
\begin{align}
\label{in1.1}
&(\triangle+k^2)u(x, k)=\mathcal{S}(k, x), x\in \mathbb{R}^d,\\
&\lim_{r\rightarrow \infty}r^{\frac{d-1}{2}}{\big\{}\frac{\partial u}{\partial r}-iku{\big\}}=0, \,\, r=|x|,\label{in1.2}
\end{align}
where $k=\omega/c_0$ is the wave number, $\omega$ is the radial frequency,  $c_0$ is the speed of sound, $d$ is the spatial dimension and \eqref{in1.2} is the Sommerfeld radiation condition.
Here, we discuss the case of $d=2$ and the source function $\mathcal{S}$ is separable with respect to space and frequency, i.e.,
\begin{align}\label{in1.3}
\mathcal{S}(k, x)=A(k)f(x),\,\, \text{supp}f\subset D.
\end{align}
In real settings, the source $f$ is usually unknown and needs to be determined from some observed data. And the data receiver devices are generally located on a remote closed surface $\partial\Omega$. We can collect the information of $u$ on $\partial\Omega$ for many frequencies $\omega$ or wave numbers $k$.  
 For $x\in\partial\Omega$, we have the potential representation of acoustic pressure  \cite{eller}
\begin{align}
u(k, x)&=\int_D \mathcal{S}(k, y)H_0^{(1)}(k|x-y|)dy\nonumber\\
&=A(k)\int_D f(y)H_0^{(1)}(k|x-y|)dy,\label{in1.4}
\end{align}
where $H_0^{(1)}$ is the cylindrical Hankel function. The data are collected on finite points $\vec{x}:=\{x_1, x_2, \cdots, x_N\}\subset\partial\Omega$ for finite wave numbers $\{k_1, k_2, \cdots, k_M\}\subset [k_{\min}, k_{\max}]$. In Figure \ref{fig:digit}, we show the measure points by the circles on the boundary of the square. By virtue of \eqref{in1.4}, we can write
\begin{align}\label{in1.5}
\mathcal{H}_{k_m}f+\eta_m=[u(x_j, k_m)]_{j=1}^N:=b_m \,\, \text{for}\,\, m=1, 2, \cdots M,
\end{align}
where $\mathcal{H}_{k}: L^2(D)\rightarrow \mathbb{C}^{N}$  depends on wave number $k$ and $\eta_m$ is the noise. Denote $\mathcal{H}:=[\mathcal{H}_{k}]_{k=k_1}^{k_M}$ and $b=[b_1; b_2; \cdots; b_M]$. 
Further assume that the function $f$ is known a priori to have the form
\begin{align}
\label{lev2.1}
f(x)=\sum_{l=1}^n w_l \mathbb{I}_{D_l}(x);
\end{align}
here $\mathbb{I}_D$ denotes the indicator function of subset $D\subset\mathbb{R}^2$, $\{D_l\}_{l=1}^n$ are subsets of $D$ such that
$\bigcup_{l=1}^n\bar{D}_l=\bar{D}$ and $D_l\bigcap D_j=\varnothing$ ($l\neq j$), the $\{w_l\}_{l=1}^n$ are known positive constants. In this setting, the regions $D_l$ determine the unknown function and therefore become the primary unknowns. 
This problem has been discussed in \cite{alzaalig,eller}, from which we can see that the source $f$ can be determined  uniquely according to the multiple frequency data. Many very successful approaches have been proposed to solve the reconstruction problem of the acoustic source \cite{bao, griesmaier_1,  griesmaier_2, liu, potthast,  sun, wang_2, zhang}. The main techniques in these works include iteration and sampling algorithms. And coping with the ill-posedness is usually addressed through the use of regularization. 

 For this kind of geometry shape problems, the other very powerful tool is the level set method \cite{osher1}, which was originally introduced by Osher and Sethian \cite{osher2}. It is devised as a simple and versatile method for computing and analyzing the motion of an interface. 
In this method, the interface is supposed to move in the normal direction with some speed, and then a partial differential equation, in particular a Hamilton-Jacobi equation, about the level set function is established.
Due to the high efficiency, it has developed to be one of the most popular tools for many disciplines, such as image processing, computer graphics, computational geometry. Besides the aspect of high efficiency in computational schemes, the crucial advantage of the level set approach lies in that it allows for topological changes to be detected during the course of algorithms, which is impossible with classical methods based on curve parameterizations. For shape reconstruction inverse problems, level set method  has also been widely studied \cite{aghasi, burger, dorn, ito, santosa1}.

Recently, statistical perspective has been widely applied in the field of inverse problems \cite{kaipio_1, stuart, tarantola}. In this framework, the unknown, noise and data are viewed as random variables. Moreover, the core issue for inverse problems in the statistics framework is how to characterize the posterior distribution, e.g., evaluate posterior moments or other posterior expectations. Based on the Bayes' formula, plenty of algorithms have been proposed to characterize the posterior distribution. 
  The most widely used class of approaches include Markov chain Monte Carlo (MCMC) \cite{stuart, wang_1}, maximum a posterior estimation (MAP) \cite{stuart}, ensemble Kalman filter (EnKF) \cite{iglesias2} etc. These methods has been applied in many inverse problems, e.g.,  the electrical impedance tomography \cite{dunlop_2, kaipio_2, kaipio_3, roininen}, the Darcy flow \cite{beskos, iglesias4}, fluid mechanics \cite{cotter, iglesias2}, optical flow \cite{bereziat},  inverse scattering  \cite{bui-thanh} etc. A natural idea to solve the interface reconstruction problem is to combine the statistical perspective and the level set approach. To the authors' knowledge, the studies in this field are rare. 
Since the level set mapping is not continuous, it is difficult to establish the well-posedness of the Bayesian level set posterior distribution according to the theory in \cite{stuart}.
To get the well-posedness, the Gaussian prior assumption on the level set function  is imposed in  \cite{iglesias1} and the almost surely continuous of the forward operator on the level set function is obtained. And subsequently the posterior well-posedness is given \cite{iglesias1}.  In \cite{dunlop_1}, the hierarchical level set Bayesian algorithm is proposed to solve the kind of problems with unknown prior parameter. In \cite{chada}, an EnKF coupling level set is considered and applied in EIT problem. One can refer to  \cite{lorentzen1, lorentzen2, ping, tai} for some more real applications, e.g., oil reservoir problem, image segmentation problem etc.





\section{Level set approach}

In this section, we sketch the level set approach for inverse shape reconstruction problems.


According to \cite{dorn,dunlop_1,iglesias1}, we can characterize the region $D_l$  by the following way
\begin{align}\label{lev2.2}
D_l=\{x\in D\mid c_{l-1}\leq\varphi(x)<c_{l}\},
\end{align}
where $\varphi$ is called the level set function and $-\infty=c_0<c_1<\cdots<c_n=\infty$ are constants. 
Using this characterization and the assumption \eqref{lev2.1}, we can define
the level set map $\mathcal{G}: C(\bar{D})\rightarrow L^2(D)$ 
\begin{align}\label{lev2.3}
\mathcal{G}(\varphi)(x)=f(x)=\sum_{l=1}^n w_l\mathbb{I}_{D_l}(x).
\end{align}
It is clear that the mapping $\mathcal{G}$ is not an injection since there exist many different $\varphi$ determine the same shape.  However, every given function $\varphi$ uniques specifies a function $f$ \cite{dorn}. 
We combine \eqref{lev2.3} and \eqref{in1.5} to give
\begin{align}\label{lev2.4}
\mathcal{H}_{k_m}\circ\mathcal{G}(\varphi)=[u(k_m, x_j)]_{j=1}^N \,\, \text{for}\,\, m=1, 2, \cdots M.
\end{align}
This can be rewritten in a compact form
\begin{align}
\label{lev2.5}
\mathcal{K}(\varphi)+\eta=b,
\end{align}
where $\mathcal{K}=[\mathcal{H}_{k_m}\circ\mathcal{G}]_{m=1}^M$ and $b=[u(k_m, x_j)]$ is an $N\times M$ dimensional column vector.

Introduce a time variable $t$ and define 
\begin{align}
\label{leva1}
\partial D(t)=\{x: \varphi(x, t)=0\}.
\end{align}
The level set function $\varphi$ evolves dynamically given by the convection equation \cite{burger,santosa1}
\begin{align}
\label{lev2.6}
\frac{\partial\varphi}{\partial t}=-v\cdot\nabla\varphi=-v|\nabla\varphi|\cdot\frac{\nabla\varphi}{|\nabla\varphi|}:=-\mathfrak{v}|\nabla\varphi|,
\end{align}
which reflects that $\varphi$ varies along with the direction of the negative gradient with velocity $v$. According to \cite{santosa1}, we take the variate $\mathfrak{v}$  as
\begin{align}\label{lev2.7}
\mathfrak{v}=J(f)^T(b-\mathcal{H}f),
\end{align}
where $J(f)$ is the Jacobian of $\mathcal{H}f$ at $f$. For an initial guess $\varphi_0$, we have the following initial value problem for $\varphi(x, t)$ of Hamilton-Jacobi system
\begin{align}
\label{lev2.8}
&\frac{\partial\varphi}{\partial t}=J(f)^T(\mathcal{H}f-b)|\nabla\varphi|,
\\
&\varphi(x, 0)=\varphi_0.
\label{lev2.9}
\end{align}
An explicit discretized version of the above system yields
\begin{align}
\label{lev2.10}
\varphi_{n+1}=\mathcal{M}(\varphi_n)
\end{align}
with initial guess $\varphi_0$.

\section{EnKF for the level set algorithm}
We now treat problem \eqref{lev2.5} in the frame of statistical inversion, which looks all involved quantities as random variables and explores the posterior distribution on $\varphi$ given $b$. We restate the inverse problem as an identification problem of system state. The system evolves in dynamics \eqref{lev2.10} with unknown initial guess. We assume the initial state $\varphi_0$ obeys the distribution $\mu_0=N(\bar{\varphi}, \Xi_0)$, where $\bar{\varphi}$ is the prior mean and $\Xi_0$ is the prior covariance operator. Observations $b$ are collected at time $t_n$ denoted by $b^n$. Thereby, we have the following state identification problem
\begin{align}
&\varphi_{n+1}=\mathcal{M}(\varphi_n),\,\, n=0, 1, \cdots, N_t\label{st3.1}\\
&\varphi_0\sim\mu_0,\label{st3.2}\\
&\mathcal{K}(\varphi_{n+1})+\eta^{n+1}=b^{n+1},\,\, n=0, 1, \cdots, N_t,\label{st3.3}
\end{align}
where $\{\eta^{n+1}\}$ are the i.i.d.  Gaussian noise, $\eta^1\sim N(0, \Gamma)$ and $\Gamma$ is the noise covariance matrix.
We focus on the posterior distribution $\mu_{0\mid 1:N_t}$, the distribution of $\varphi_0|b^{1:N_t}$, $\varphi_0$ given $b^1, b^2, \cdots, b^{N_t}$. Bayes rule gives a characterization of $\mu_{0\mid 1:N_t}$ via the ratio of its density with respect to that of the prior
\begin{align}\label{st3.4}
\frac{d\mu_{0\mid 1:N_t}}{d\mu_0}(\varphi_0)=\frac{d\mu(b^{1:N_t}\mid\varphi_0)}{d\mu(b^{1:N_t})}.
\end{align}
And therefore, the posterior density satisfies
\begin{align}\label{st3.5}
\pi_{0\mid 1:N_t}\propto \pi(b^{1:N_t}\mid\varphi_0)\pi_0.
\end{align}
From the angle of filter, we are concerned with obtaining the posterior distribution $\mu_{n+1|n+1}$ via $\mu_{n|n}$, the distribution associated with the probability measure on the random variable $\varphi_n| b^n$.
Suppose that we get the present system state $\varphi_n$ together with the prior covariance $\Xi_n$. The operators $\Xi_n$ characterize model uncertainty and are design parameters. With these information as the prior for the next state $\varphi_{n+1}$, we estimate the new state $\varphi_{n+1}$ by blending the observed data $b^{n+1}$ into the dynamics. However, there is, in general, no easily usable closed form expression for the density of $\mu_{n+1|n+1}$. We can break this into a two-step process: prediction step and analysis step. In the prediction step, we use the present state distribution to predict the new state according to the dynamics as
\begin{align}
\label{st3.7}
\hat{\varphi}_{n+1|n}=\mathcal{M}(\varphi_n).
\end{align}
Using the prediction $\hat{\varphi}_{n+1|n}$ together with the design parameter $\Xi_n$, we give the prior distribution $\hat{\mu}_{n+1|n}$ and  further obtain the analysis state distribution according to the observational process \eqref{st3.3}
\begin{align}
\label{st3.8}
\frac{d{\mu}_{n+1|n+1}}{d\hat{\mu}_{n+1|n}}(\varphi)=\exp(-\frac{1}{2}|\mathcal{K}(\varphi)-b^{n+1}|_{\Gamma}^2).
\end{align}
When the design parameter $\Xi_n$ is fixed, this filter corresponds to the 3DVar process. Naturally, we also can estimate the prior covariance $\Xi_n$ via samples.  To characterize the posterior distribution, the ensemble Kalman filter (EnKF) approach represents the filtering distribution through an ensemble of particles. The prediction sample can be obtained for each particle. Then the prediction sample mean is taken as the prior mean and the covariance of prediction samples is taken as the prior covariance.  
In details, suppose that at  the $n$-th step, we have samples $\{\varphi_n^{(j)}\}_{j=1}^J$, $J$ is the number of samples. For each particle, we have a prediction 
\begin{align}\label{st3.10}
\hat{\varphi}_{n+1}^{(j)}=\mathcal{M}(\varphi_n^{(j)}),\,\, j=1:J.
\end{align}
 Define the empirical covariance $\hat{\Xi}_n$ as
\begin{align}\label{st3.11}
\hat{\Xi}_n=\frac{1}{J}\sum_{j=1}^{J}(\hat{\varphi}_{n+1}^{(j)}-\bar{\varphi}_{n+1})\otimes(\hat{\varphi}_{n+1}^{(j)}-\bar{\varphi}_{n+1}),
\end{align}
where $\bar{\varphi}_{n+1}=\frac{1}{J}\sum_{j=1}^{J}\hat{\varphi}_{n+1}^{(j)}$ denotes the ensemble mean.
 The analysis state for each particle is then obtained by minimizing the following functional
\begin{align}\label{st3.12}
\mathfrak{J}(\varphi)&=\frac{1}{2}\int_{\Omega}(\varphi-\hat{\varphi}_{n+1}^{(j)})\hat{\Xi}_n^{-1}(\varphi-\hat{\varphi}_{n+1}^{(j)})dx\nonumber\\
&+\frac{1}{2}|\mathcal{K}(\varphi)-b^{n+1}|_\Gamma^2,\,\, j=1:J.
\end{align}
Due to the nonlinearity of the observational operator $\mathcal{K}$ in \eqref{st3.12},  the numerical implementation of the minimization becomes difficulty. To overcome this point,  iterative EnKF approaches are proposed in \cite{iglesias2, iglesias3, yang1}. 
Define the state vector $\varsigma=[\varphi, f]^T$.
We modify the prediction process as
\begin{align}\label{sp3.11}
\hat{\varsigma}_{n+1}=\Theta(\varsigma_n)
\end{align}
 by
 \begin{align}\label{sp3.12}
\begin{split}
&\hat{\varphi}_{n+1}=\mathcal{M}(\varphi_n),\\
&\hat{f}_{n+1}=\mathcal{G}(\hat{\varphi}_{n+1}).
\end{split}
\end{align}
Accordingly, the observation operator is given by $\mathcal{O}=[0, \mathcal{H}]$.
Formally, the minimization function about the state vector $\varsigma$ has the form of
\begin{align}\label{sp3.13}
\mathfrak{J}(\varsigma)=\frac{1}{2}\|\varsigma-\hat{\varsigma}_{n+1}\|^2_{\hat{\Xi}_n^\varsigma}+\frac{1}{2}|\mathcal{O}\varsigma-b^{n+1}|_\Gamma^2,
\end{align}
where $\hat{\Xi}_n^\varsigma$ is the prior covariance about $\varsigma$. 
The minimizer can be given by
\begin{align}
\nonumber\varsigma&=((\hat{\Xi}_n^\varsigma)^{-1}+\mathcal{O}^*\Gamma^{-1}\mathcal{O})^{-1}((\hat{\Xi}_n^\varsigma)^{-1}
\hat{\varsigma}_{n+1}+\mathcal{O}^*\Gamma^{-1}b^{n+1})\\
\label{up3.14}
&=(I+\hat{\Xi}_n^\varsigma\mathcal{O}^*\Gamma^{-1}\mathcal{O})^{-1}(\hat{\varsigma}_{n+1}+\hat{\Xi}_n^\varsigma\mathcal{O}^*\Gamma^{-1}b^{n+1}).
\end{align}
Write the prior covariance matrix $\hat{\Xi}_n^\varsigma$ as 
\begin{align}
\hat{\Xi}_n^\varsigma=
\left(\begin{array}{cc}
 \hat{\Xi}_n^{\varphi\varphi}&   \hat{\Xi}_n^{\varphi f}
\\ \hat{\Xi}_n^{f\varphi} &  \hat{\Xi}_n^{ff}\end{array}
\right).
\end{align}
The component form of \eqref{up3.14} has
\begin{align}
\label{dyd3.16}&\varphi=\hat{\varphi}_{n+1}+\hat{\Xi}_n^{\phi f}\mathcal{H}^*(\mathcal{H}\hat{\Xi}_n^{ff}\mathcal{H}^*+\Gamma)^{-1}(b^{n+1}-\mathcal{H}\hat{f}_{n+1}),\\
\label{dyd3.17}
&f=(I+\hat{\Xi}_n^{ff}\mathcal{H}^*\Gamma^{-1}\mathcal{H})^{-1}(\hat{f}_{n+1}+\hat{\Xi}_n^{ff}\mathcal{H}^*\Gamma^{-1}b^{n+1}).
\end{align}
We list this algorithm here.
\begin{center}
Algorithm 1:
\end{center}
\begin{itemize}
\item[\labelitemi] Initial particles: Let $\{\varphi^{(j)}_0\}_{j=1}^J$ be the initial samples drawn from the prior distribution $\mathcal{N}(\bar{\varphi}, \Xi_0)$. And obtain the initial samples of $f^{(j)}_0$ according to $f^{(j)}_0=\mathcal{G}(\varphi^{(j)}_0)$.
\item[\labelitemi] Prediction: Evaluate 
\begin{align}\label{st3.13}
&\hat{\varphi}^{(j)}_n=\mathcal{M}(\varphi^{(j)}_{n-1}),\,\, j=1:J,\\
&\hat{f}^{(j)}_n=\mathcal{G}(\hat{\varphi}^{(j)}_n),\,\, j=1:J. \label{st3.14}
\end{align}
Compute $\bar{\varphi}_n=\frac{1}{J}\sum_{j=1}^J\hat{\varphi}^{(j)}_n$,
$\bar{f}_n=\frac{1}{J}\sum_{j=1}^J\hat{f}^{(j)}_n$.
\item[\labelitemi] Analysis: Define $\hat{\Xi}_n^{\varphi f}$, $\hat{\Xi}_n^{ff}$ by
\begin{align}
&\hat{\Xi}_n^{\varphi f}=\frac{1}{J}\sum_{j=1}^{J}(\hat{\varphi}^{(j)}_n-\bar{\varphi}_n)\otimes(\hat{f}^{(j)}_n-\bar{f}_n), \label{st3.15}\\
&\hat{\Xi}_n^{ff}=\frac{1}{J}\sum_{j=1}^{J}(\hat{f}^{(j)}_n-\bar{f}_n)\otimes(\hat{f}^{(j)}_n-\bar{f}_n).\label{st3.16}
\end{align}
Update each particle via \eqref{dyd3.16} and \eqref{dyd3.17}.
Compute the mean of the parameter update
\begin{align}\label{sth3.18}
\varphi_n=\frac{1}{J}\sum_{j=1}^J \varphi^{(j)}_n,\\
\label{sth3.19}
f_n=\frac{1}{J}\sum_{j=1}^J f^{(j)}_n
\end{align}
and moreover check the stopping rule.
\end{itemize}
In other way, the state vector is defined by $\zeta:=[\varphi, \mathfrak{b}]^T$ \cite{iglesias2}. Accordingly, the observation operator is given as $[0, I]$.
In this setting, the EnKF approach is listed in the follows. 

\begin{center}
Algorithm 2:
\end{center}
\begin{enumerate}
\item[\labelitemiii] Initial particles: Let $\{\varphi^{(j)}_0\}_{j=1}^J$ be the initial samples drawn from the prior distribution $\mathcal{N}(\bar{\varphi}, \Xi_0)$. 
\item[\labelitemiii] Prediction: Evaluate 
\begin{align}\label{st3.24}
&\hat{\varphi}^{(j)}_n=\mathcal{M}(\varphi^{(j)}_{n-1}),\,\, j=1:J,\\
&\hat{\mathfrak{b}}_{(j)}^n=\mathcal{K}(\hat{\varphi}^{(j)}_n).
\end{align}
Compute $\bar{\varphi}_n=\frac{1}{J}\sum_{j=1}^J\hat{\varphi}^{(j)}_n$,
$\bar{\mathfrak{b}}^n=\frac{1}{J}\sum_{j=1}^J\hat{\mathfrak{b}}_{(j)}^n$.
\item[\labelitemiii] Analysis: Define $\hat{\Xi}_n^{\mathfrak{b}}$, $\hat{\Xi}_n^{\mathfrak{b}\mathfrak{b}}$ by
\begin{align}
&\hat{\Xi}_n^{\varphi\mathfrak{b}}=\frac{1}{J}\sum_{j=1}^{J}(\hat{\varphi}^{(j)}_n-\bar{\varphi}_n)\otimes(\hat{\mathfrak{b}}_{(j)}^n-\bar{\mathfrak{b}}^n)\label{st3.15}\\
&\hat{\Xi}_n^{\mathfrak{b}\mathfrak{b}}=\frac{1}{J}\sum_{j=1}^{J}(\hat{\mathfrak{b}}_{(j)}^n-\bar{\mathfrak{b}}^n)\otimes(\hat{\mathfrak{b}}_{(j)}^n-\bar{\mathfrak{b}}^n).
\end{align}
Update each particle by 
\begin{align}\label{st3.17}
\varphi^{(j)}_n=\hat{\varphi}^{(j)}_n+\hat{\Xi}_n^{\varphi\mathfrak{b}}(\hat{\Xi}_n^{\mathfrak{b}\mathfrak{b}}+\Gamma)^{-1}(b^n-\hat{\mathfrak{b}}_{(j)}^n).
\end{align}
Compute the mean of the parameter update
\begin{align}\label{st3.18}
\varphi_n=\frac{1}{J}\sum_{j=1}^J \varphi^{(j)}_n.
\end{align}
and moreover check the stopping rule.
\end{enumerate}




\section{Prior distribution}

From the above description, we see that the prior distribution $\mu_0=N(\bar{\varphi}, \Xi_0)$ plays key role in the statistical inversion. 
In this paper, we choose the Whittle-Mat\'ern fields as the prior distribution. They are stationary Gaussian random fields with the autocorrelation function
\begin{align}
\label{pr4.1}
C(x)=\frac{2^{1-\nu}}{\Gamma(\nu)}(\frac{|x|}{l})^\nu K_\nu(\frac{|x|}{l}), \,\, x\in\mathbb{R}^d,
\end{align}
where $\nu>0$ is the smoothness parameter, $K_\nu$ is the modified Bessel function of the second kind of order $\nu$ and $l$ is the length-scale parameter. The integer value of $\nu$ determines the mean square differentiability of the underlying process, which matters for predictions made using such a model.
This distribution is discussed by many authors \cite{lindgren, rasmussen, roininen}.  A method for explicit, and computationally efficient, continuous Markov representations of Gaussian Mat\'ern fields is derived  \cite{lindgren}.
The method is based on the fact that a Gaussian Mat\'ern field on $\mathbb{R}^d$
can be viewed as a solution to the stochastic partial differential equation (SPDE)
\begin{align}\label{pr4.2}
(I-l^2\triangle)^{(\nu+d/2)/2}\varphi=\sqrt{\alpha l^2}W,
\end{align}
where $W$ is the Gaussian white noise and the constant $\alpha$ is 
\begin{align*}
\alpha:=\sigma^2\frac{2^d\pi^{d/2}\Gamma(\nu+d/2)}{\Gamma(\nu)}.
\end{align*}
Here $\sigma^2$ is the variance of the stationary field.
The operator $(I-l^2\triangle)^{(\nu+1)/2}$ is a pseudo-differential operator defined by its Fourier transform. When $\nu\in\mathbb{Z}$, we can solve \eqref{pr4.2} by the finite element method with suitable boundary condition. The stochastic weak solution of the SPDE \eqref{pr4.2} is found by
requiring that
\begin{align}
(\psi, (I-l^2\triangle)^{(\nu+d/2)/2}\varphi)=(\psi, \sqrt{\alpha l^2}W),
\end{align}
where $\psi\in L^2(D_1)$. We use the linear element to solve the weak formulation with Dirichlet boundary condition, i.e., $\varphi\mid_{\partial D_1}=0$. 

\section{Numerical test}
\begin{figure}[htbp]
\centering
\includegraphics[width=0.4\textwidth]{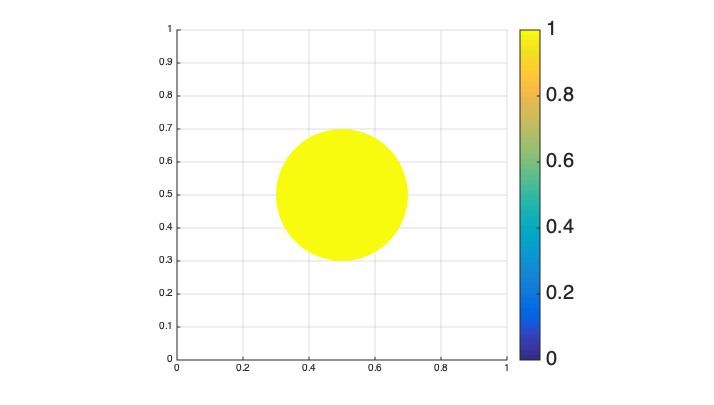}
\includegraphics[width=0.4\textwidth]{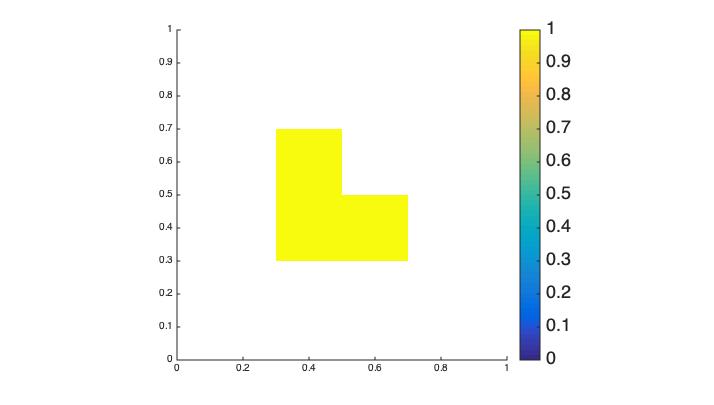}\\
\includegraphics[width=0.4\textwidth]{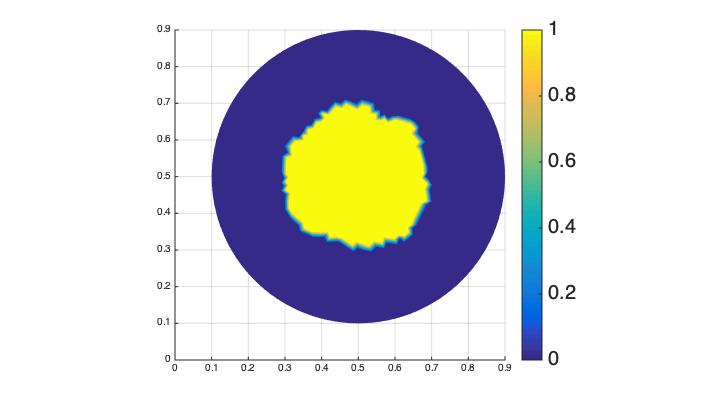}
\includegraphics[width=0.4\textwidth]{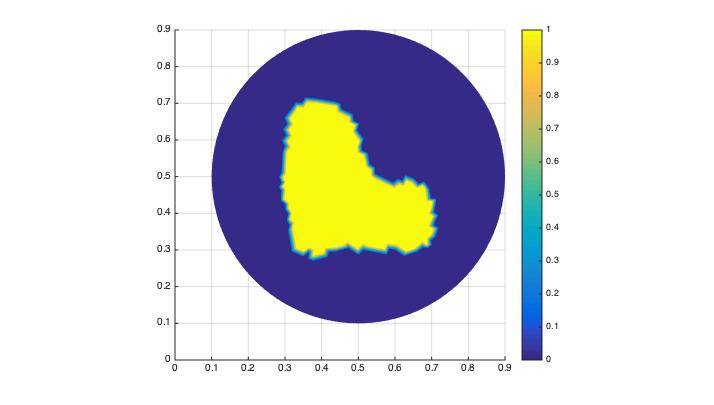}\\
\includegraphics[width=0.4\textwidth]{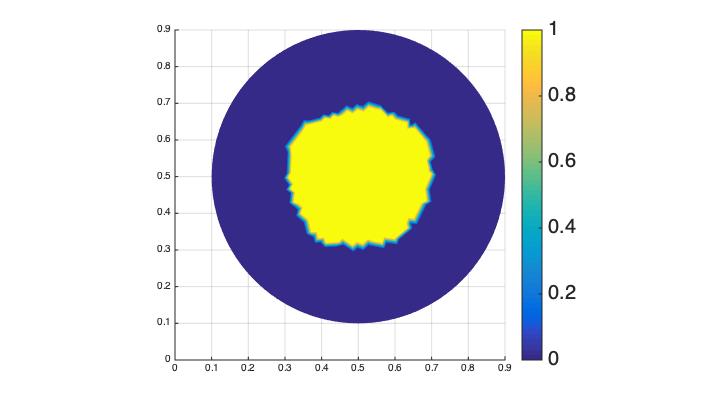}
\includegraphics[width=0.4\textwidth]{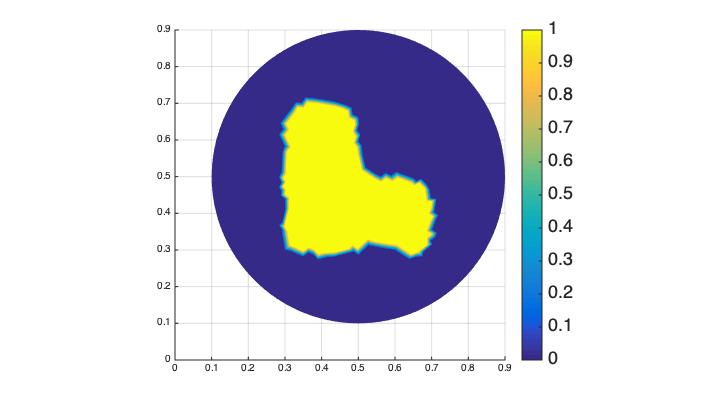}\\
\includegraphics[width=0.4\textwidth]{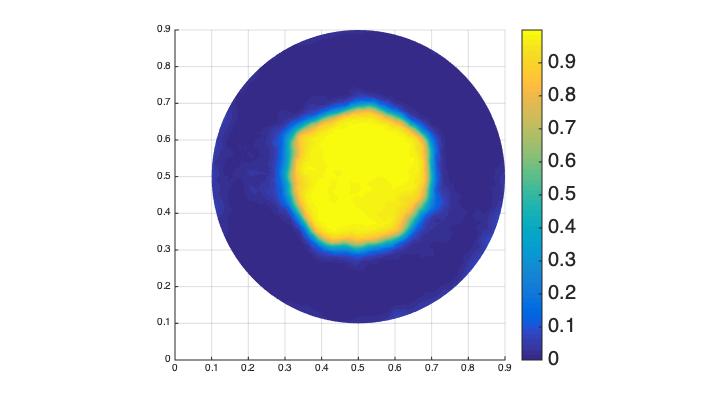}
\includegraphics[width=0.4\textwidth]{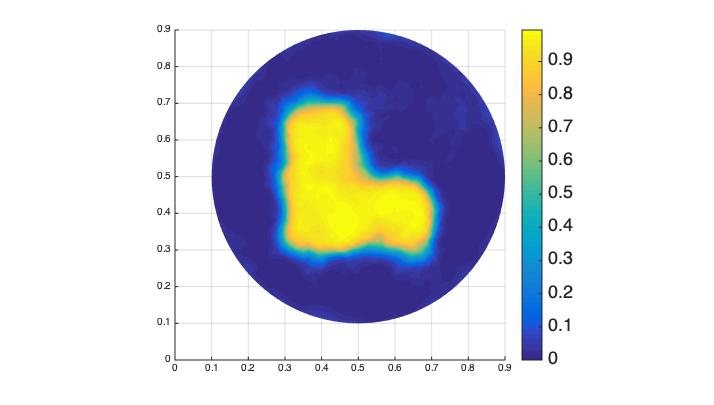}\\
\includegraphics[width=0.4\textwidth]{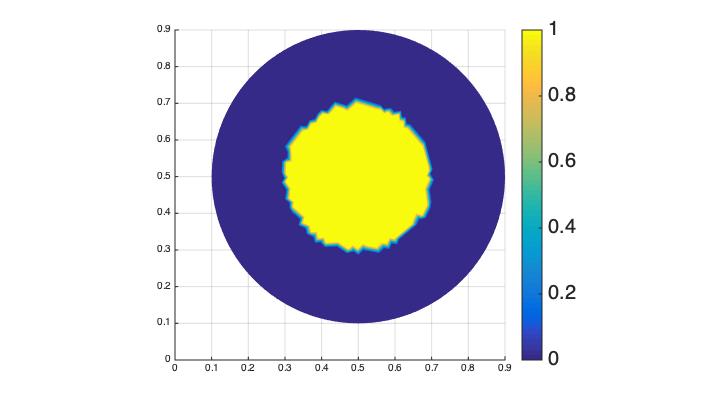}
\includegraphics[width=0.4\textwidth]{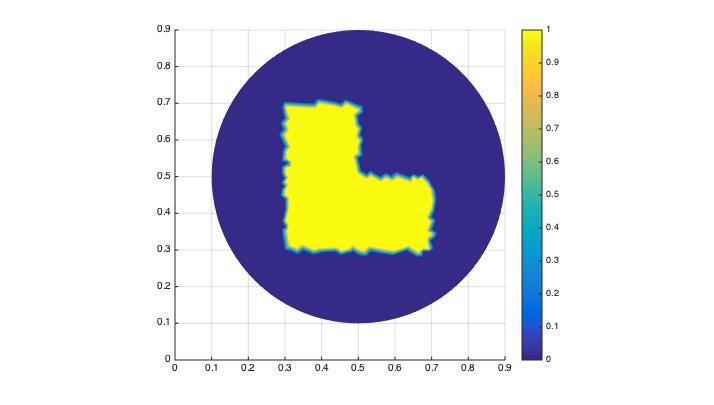}\\
\includegraphics[width=0.4\textwidth]{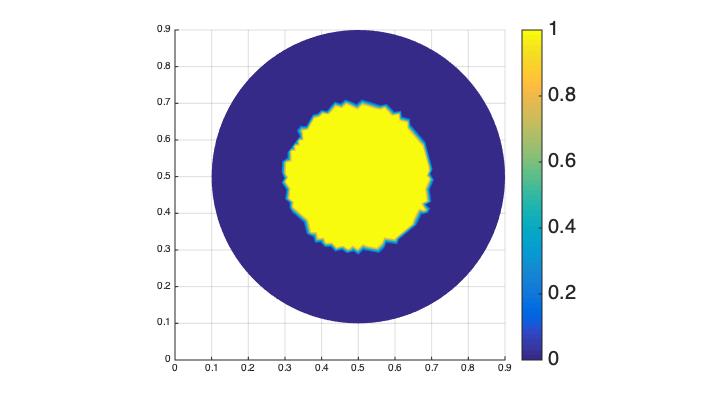}
\includegraphics[width=0.4\textwidth]{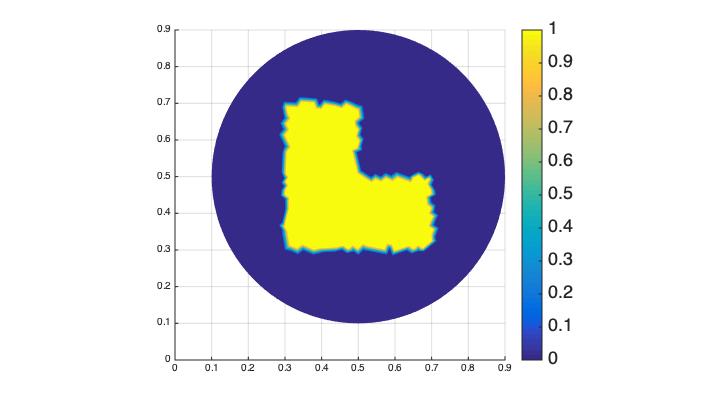}\\
\includegraphics[width=0.4\textwidth]{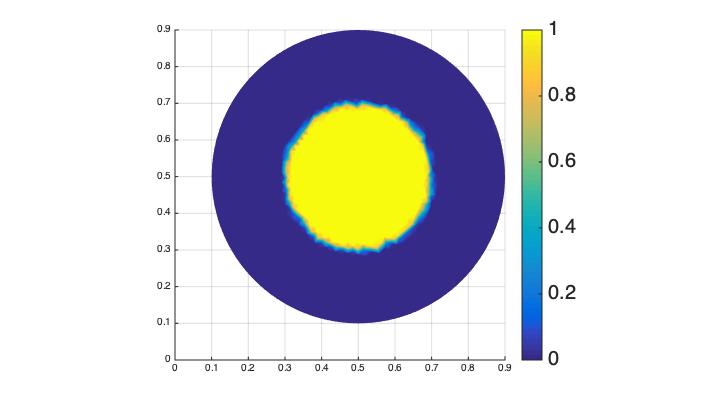}
\includegraphics[width=0.4\textwidth]{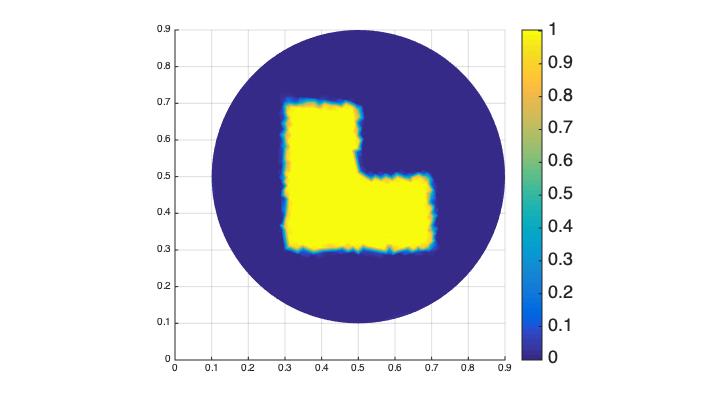}
\caption{Single domain: Line 1: True solution $f$; Lines 2-4: The reconstructions of $f$ for $\delta=0.1$; Lines 5-7: The reconstructions of $f$ for $\delta=0.01$. Lines 2 and 5 give the approximation using Algorithm 2. Lines 3 and 6 show the reconstruction using Algorithm 1 via \eqref{sth3.18}. Lines 4 and 7 are the numerical results by \eqref{sth3.19}.
}\label{fig:2}
\end{figure}

\begin{figure}[htbp]
\centering
\includegraphics[width=0.2\textwidth]{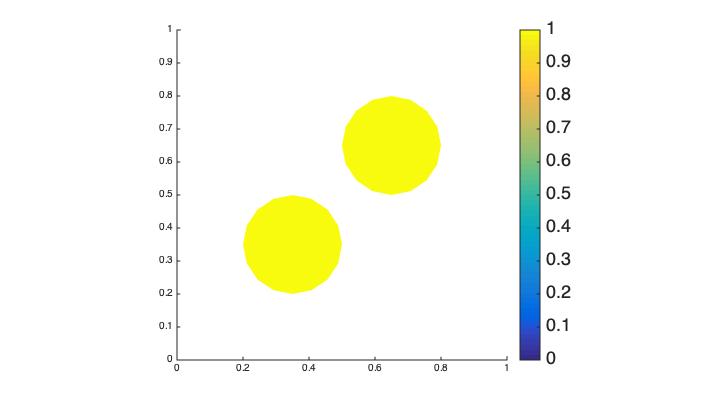}
\includegraphics[width=0.2\textwidth]{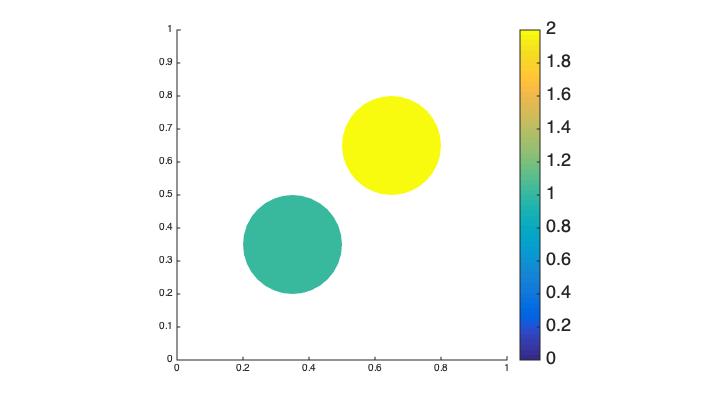}
\includegraphics[width=0.2\textwidth]{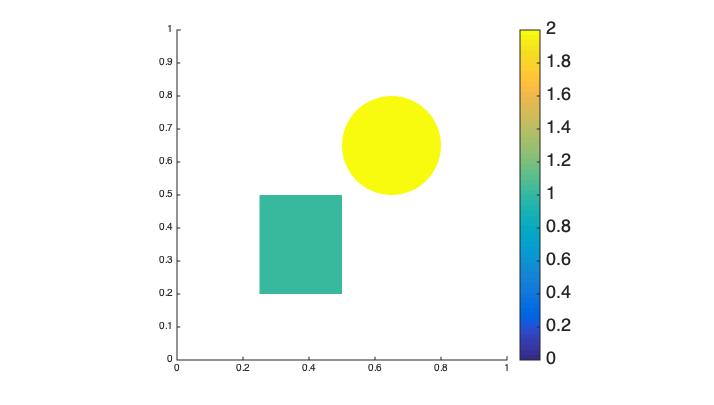}
\\
\includegraphics[width=0.2\textwidth]{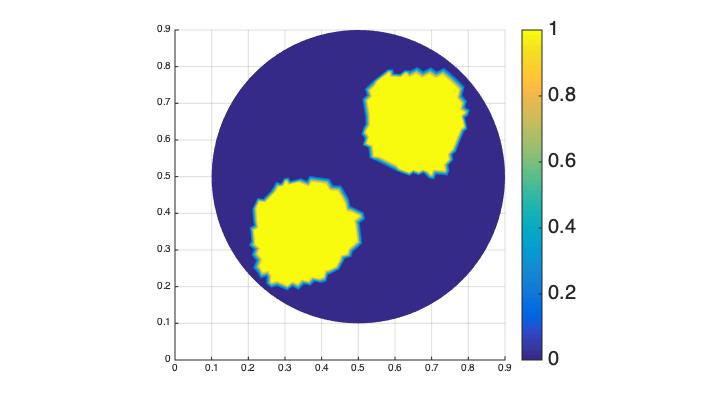}
\includegraphics[width=0.2\textwidth]{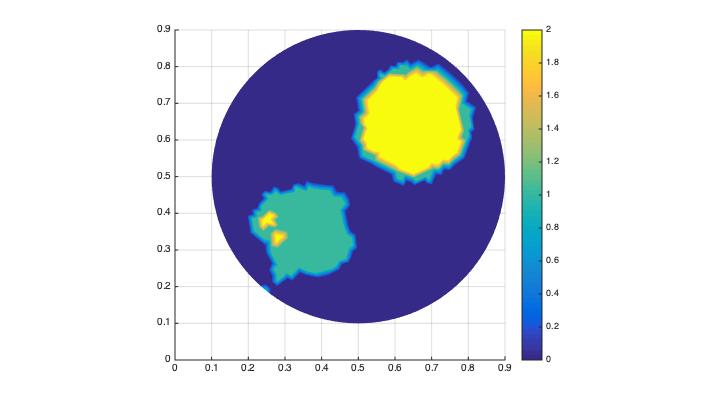}\includegraphics[width=0.2\textwidth]{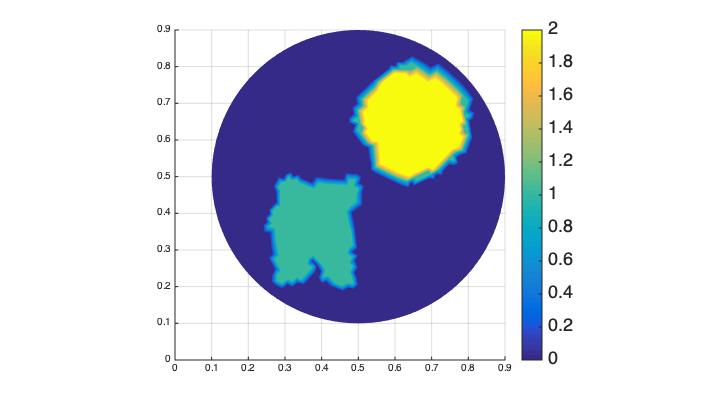}\\
\includegraphics[width=0.2\textwidth]{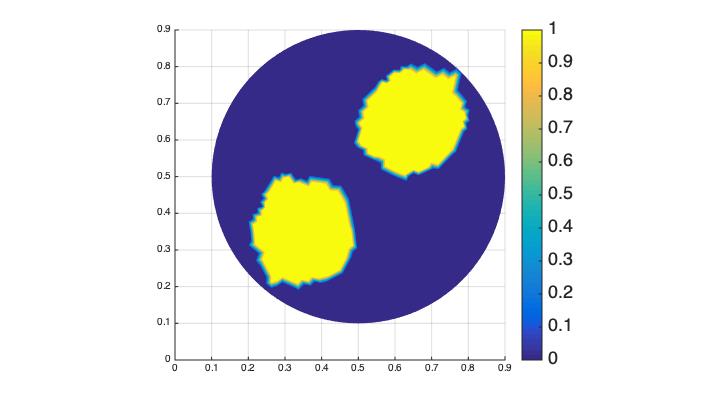}
\includegraphics[width=0.2\textwidth]{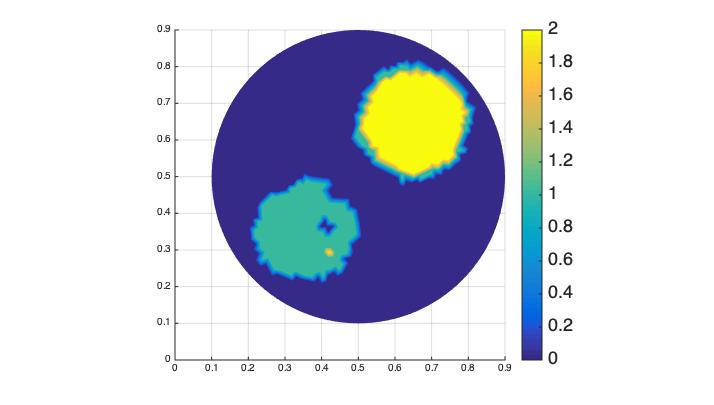}
\includegraphics[width=0.2\textwidth]{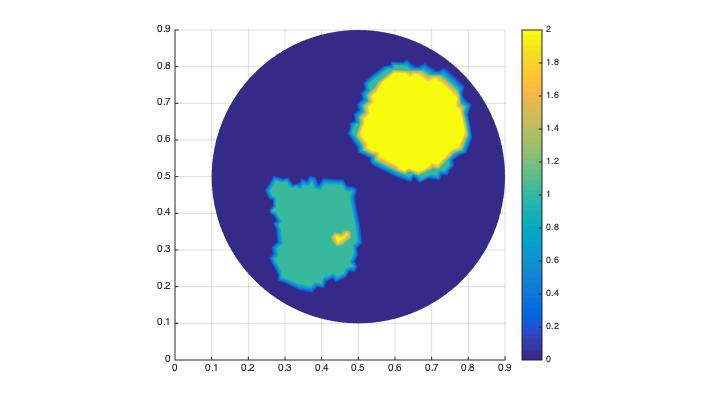}\\
\includegraphics[width=0.2\textwidth]{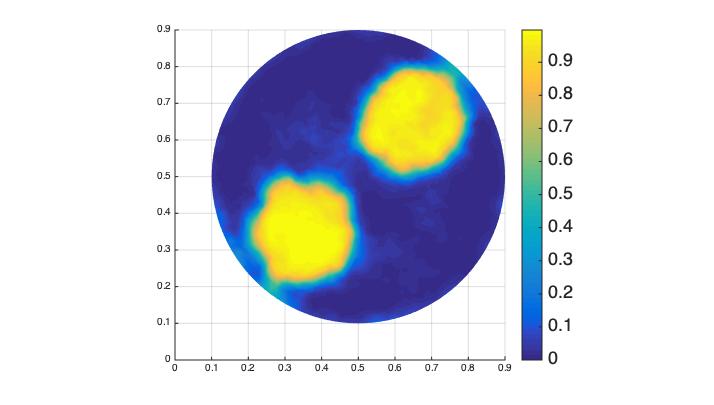}
\includegraphics[width=0.2\textwidth]{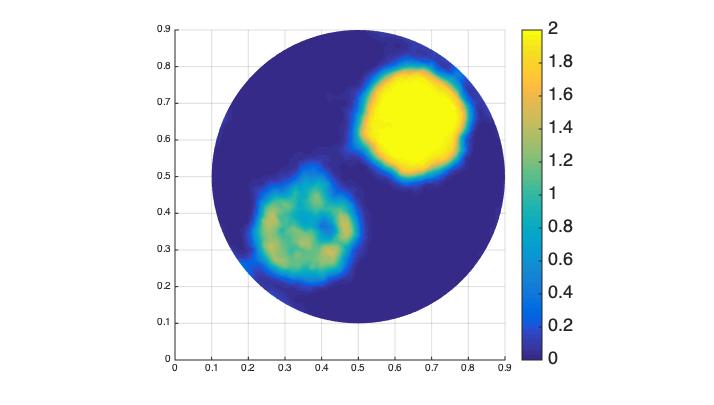}
\includegraphics[width=0.2\textwidth]{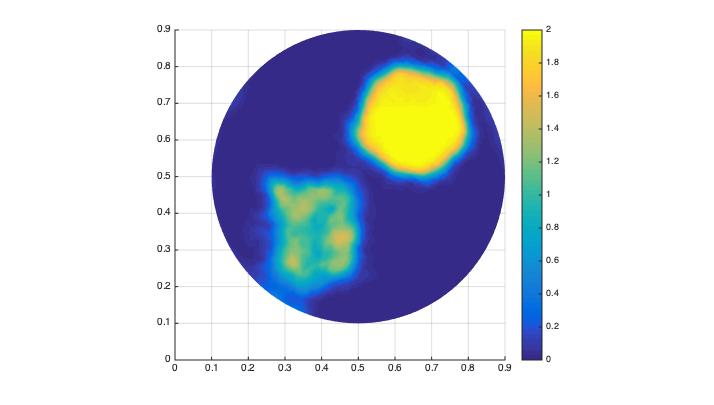}\\
\includegraphics[width=0.2\textwidth]{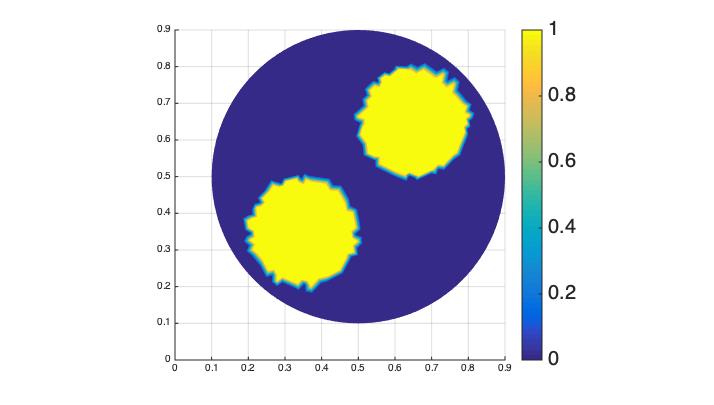}
\includegraphics[width=0.2\textwidth]{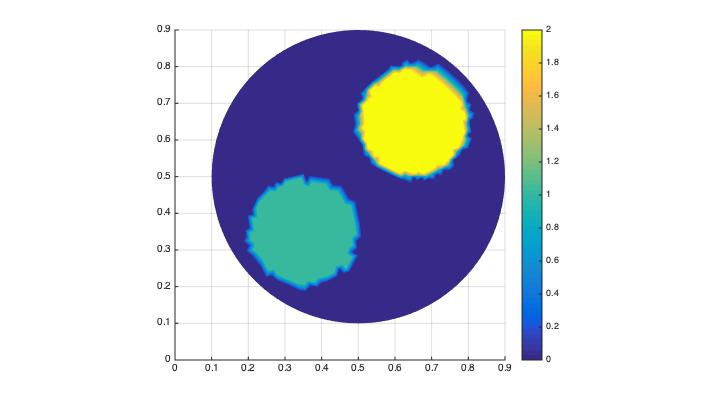}
\includegraphics[width=0.2\textwidth]{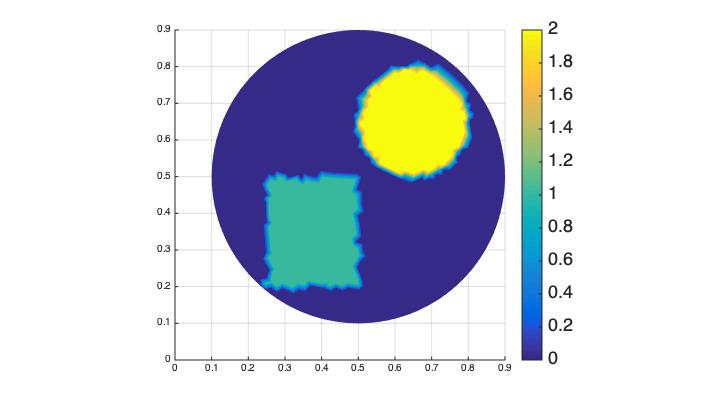}\\
\includegraphics[width=0.2\textwidth]{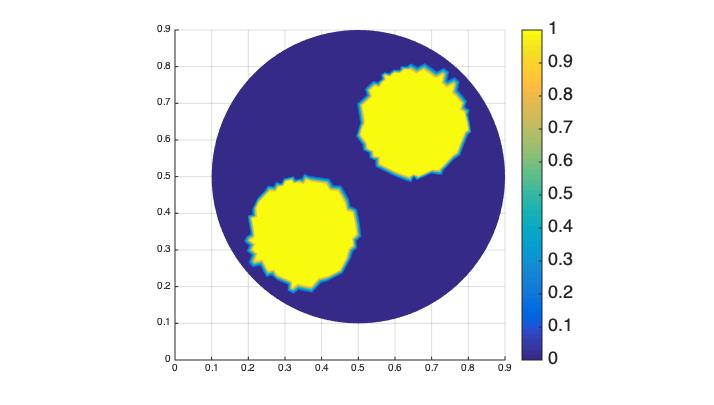}
\includegraphics[width=0.2\textwidth]{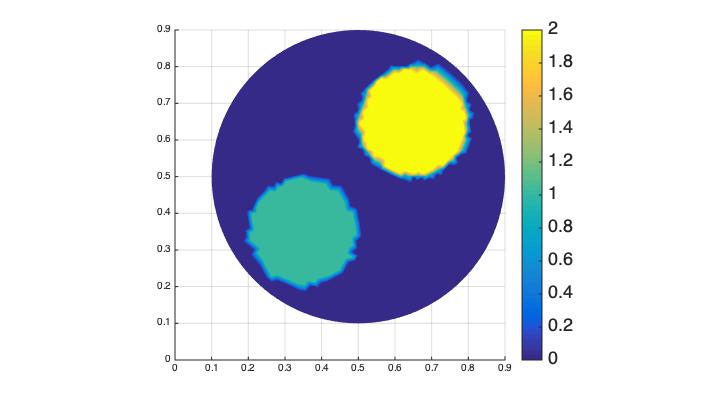}
\includegraphics[width=0.2\textwidth]{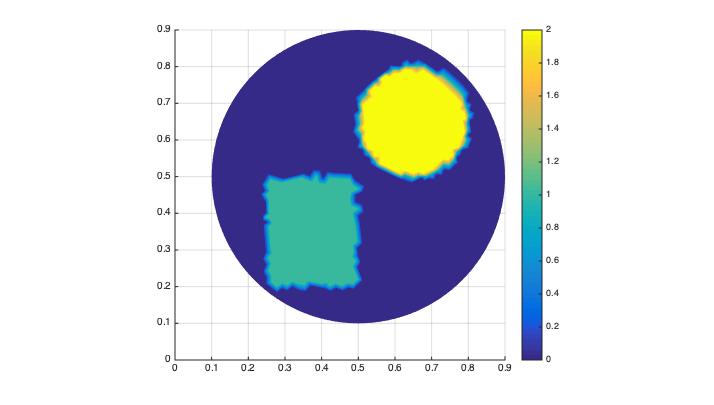}
\\
\includegraphics[width=0.2\textwidth]{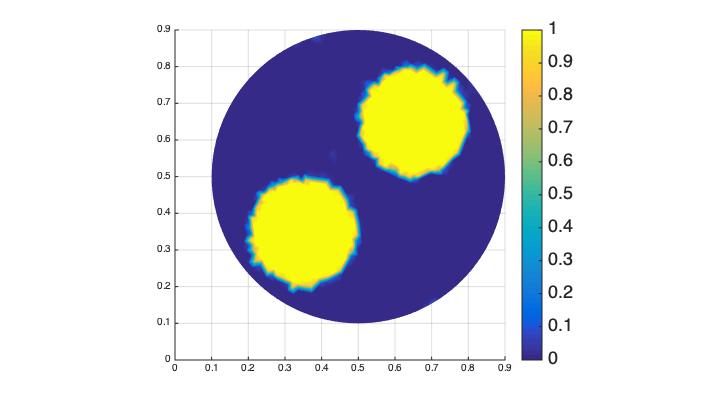}
\includegraphics[width=0.2\textwidth]{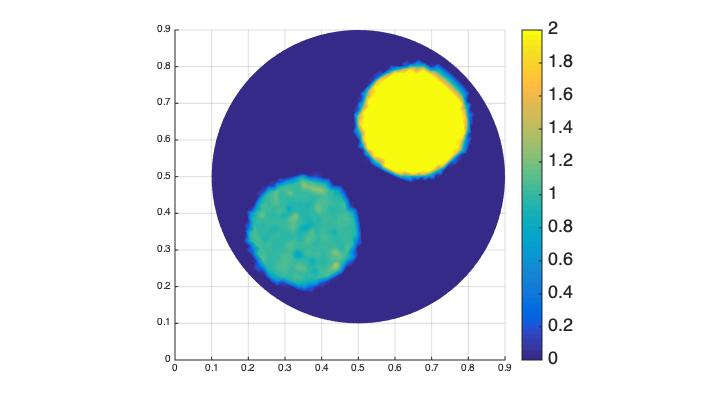}
\includegraphics[width=0.2\textwidth]{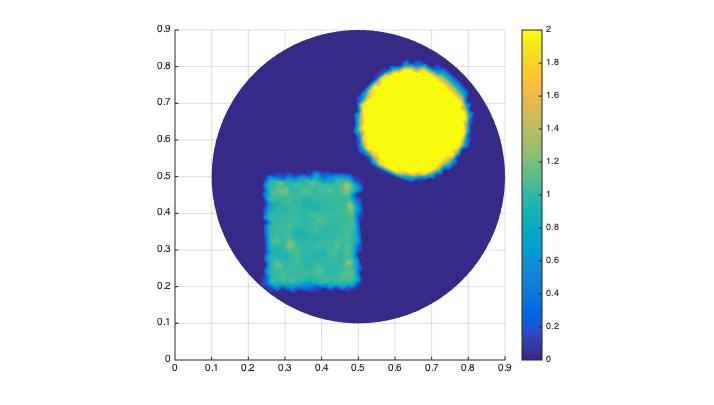}
\caption{Separate domains: Line 1: True solution $f$; Lines 2-4: The reconstructions of $f$ for $\delta=0.1$; Lines 5-7: The reconstructions of $f$ for $\delta=0.01$. Lines 2 and 5 give the approximation using Algorithm 2. Lines 3 and 6 show the reconstruction using Algorithm 1 via \eqref{sth3.18}. Lines 4 and 7 are the numerical results by \eqref{sth3.19}.}\label{fig:3}
\end{figure}

\begin{figure}[htbp]
\centering
\includegraphics[width=0.4\textwidth]{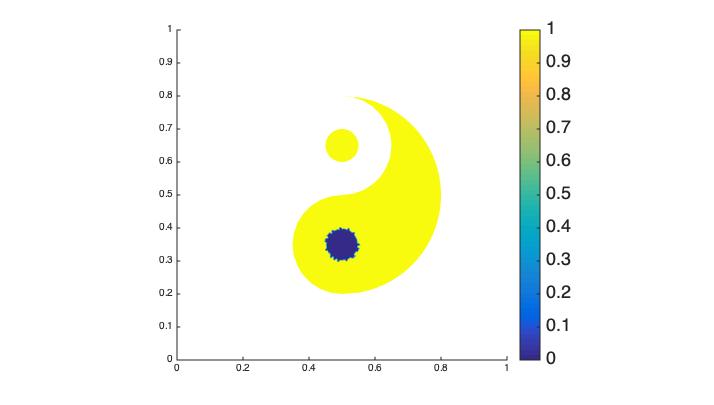}
\includegraphics[width=0.4\textwidth]{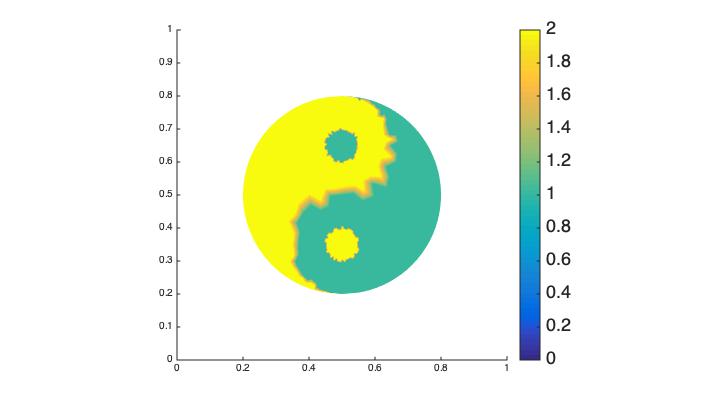}\\
\includegraphics[width=0.4\textwidth]{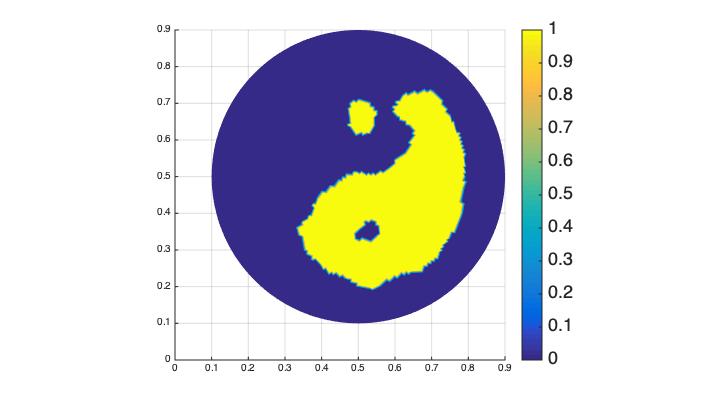}
\includegraphics[width=0.4\textwidth]{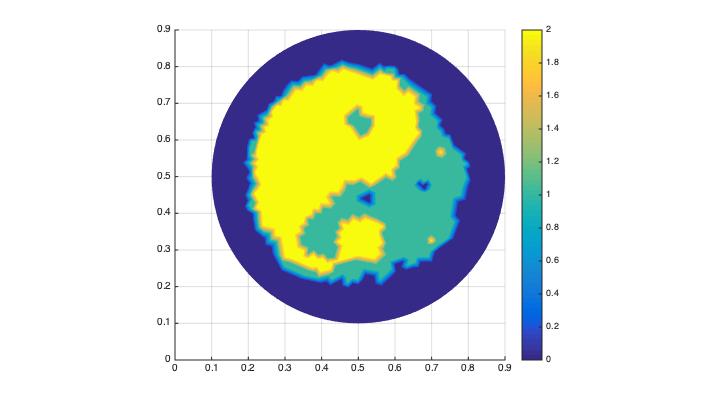}
\includegraphics[width=0.4\textwidth]{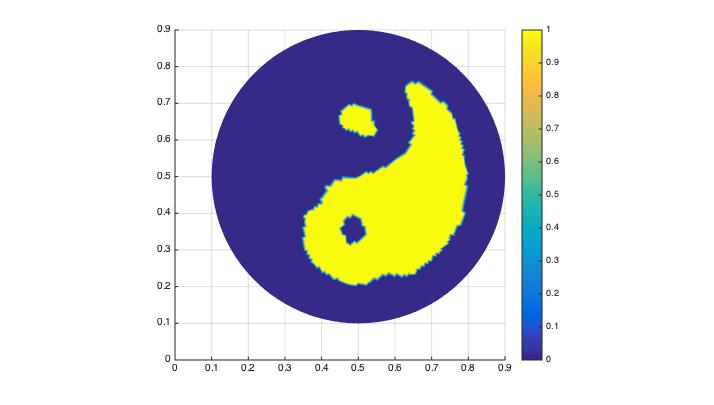}
\includegraphics[width=0.4\textwidth]{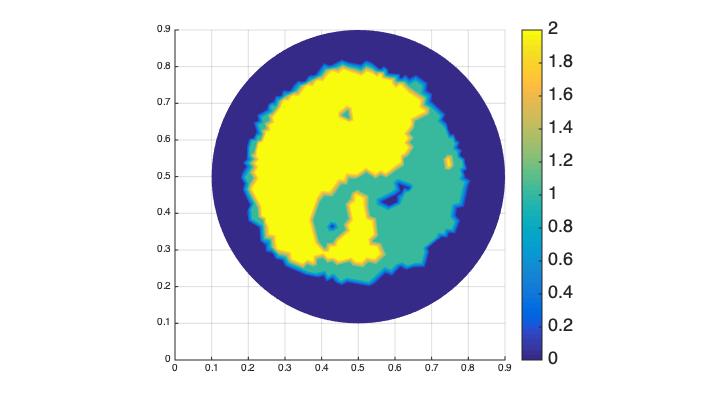}
\includegraphics[width=0.4\textwidth]{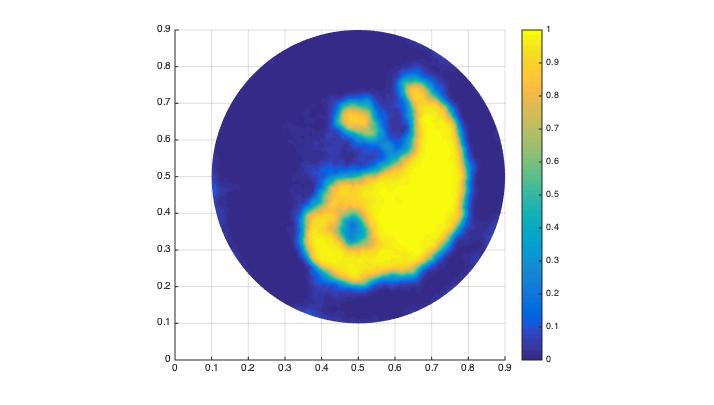}
\includegraphics[width=0.4\textwidth]{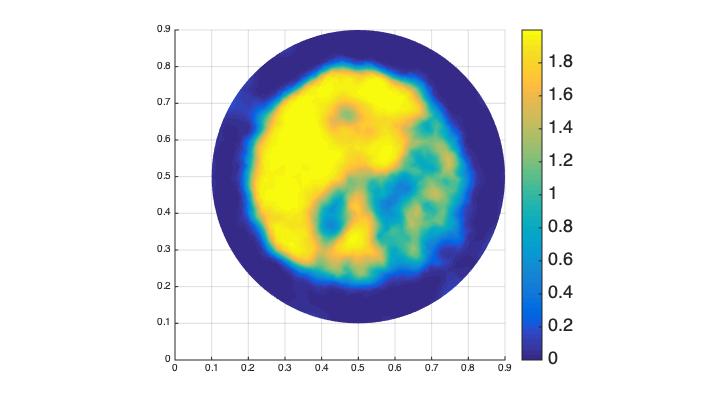}
\includegraphics[width=0.4\textwidth]{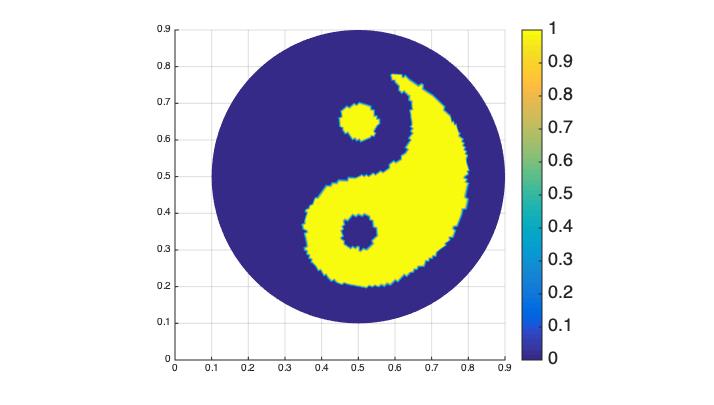}
\includegraphics[width=0.4\textwidth]{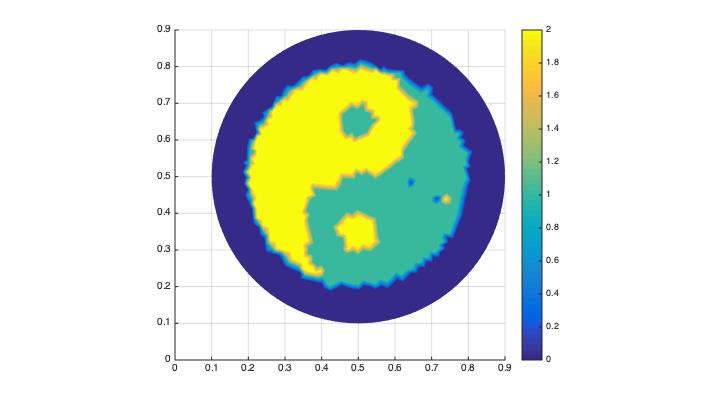}
\includegraphics[width=0.4\textwidth]{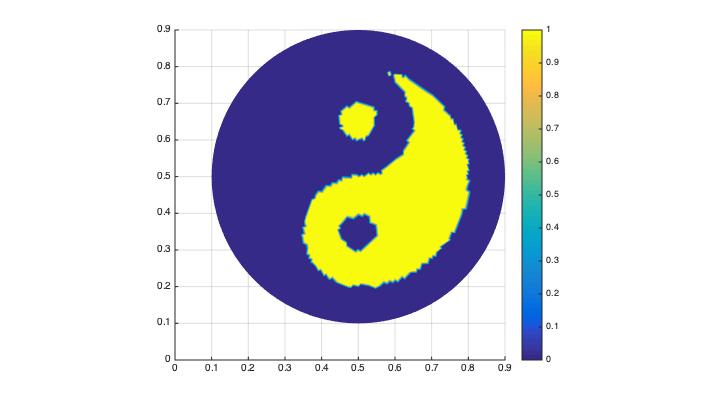}
\includegraphics[width=0.4\textwidth]{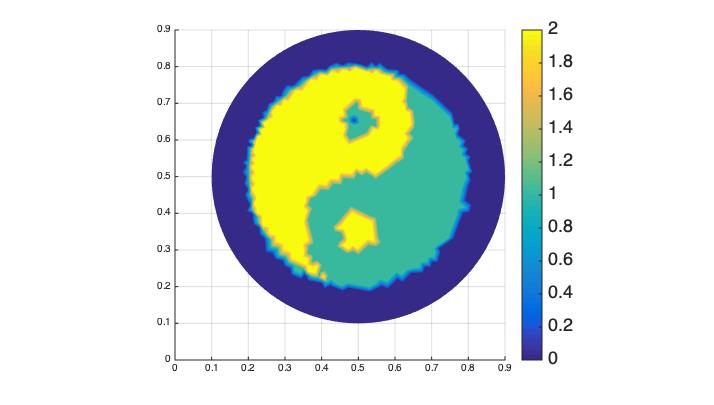}
\includegraphics[width=0.4\textwidth]{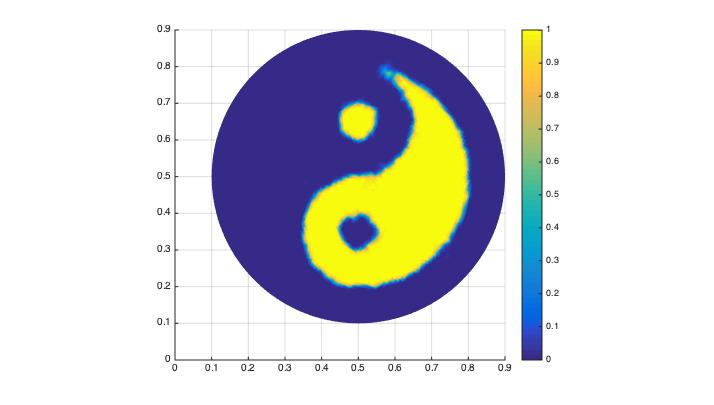}
\includegraphics[width=0.4\textwidth]{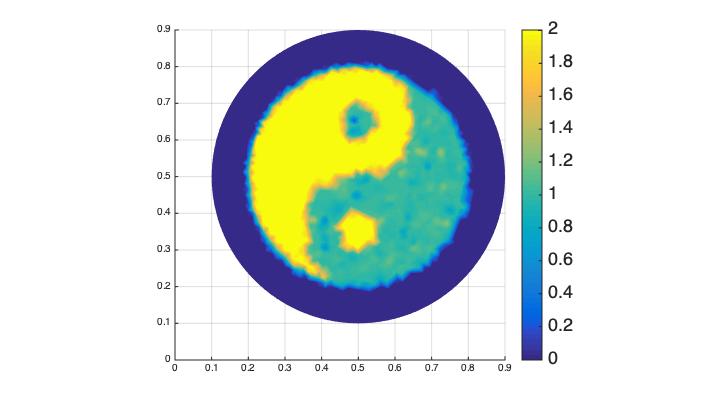}
\caption{Taichi graph: Line 1: True solution $f$; Lines 2-4: The reconstructions of $f$ for $\delta=0.1$; Lines 5-7: The reconstructions of $f$ for $\delta=0.01$. Lines 2 and 5 give the approximation using Algorithm 2. Lines 3 and 6 show the reconstruction using Algorithm 1 via \eqref{sth3.18}. Lines 4 and 7 are the numerical results by \eqref{sth3.19}.}\label{fig:4}
\end{figure}

In the numerical tests, we consider frequencies varying from $k_{\min}=50$ Hz to $k_{\max}=10$ kHz and $c_0=343$ms$^{-1}$. These parameter settings are given in \cite{eller}.  Figure \ref{fig:digit} displays the settings of our problem. The data are taken on the circle dots of the boundary $\partial\Omega$. The noise $\eta$ is taken as $\delta*\text{randn}(0, I)$. 
The prior samples are generated in the disk domain $D_1$. We assume that the support $D$ lies in the disk. The disk domain is decomposed into triangular elements and the potential integral is calculated using the Gaussian integral formulation on each triangle. 
Moreover we solve \eqref{pr4.2} in the same mesh using the linear finite element. To avoid inverse crime as far as possible, we use the finer mesh in the inverse process than in the data generation.
 The time derivative in the Hamilton-Jacobi system is discretized by the forward Euler scheme.  The Euler step size is taken as $0.1$. This sample size $J$ and the iteration step are taken as $1000$ and $100$ respectively.  The prior parameter $\nu$ is set to $1$. We display the numerical reconstruction of $f$ for several different cases. In Figures \ref{fig:2}, \ref{fig:3} and \ref{fig:4}, we give the numerical results for the cases of single domain, two separately domains and Taichi graph respectively. For each case, the data noise level are taken as $0.1$ and $0.01$. It should be noted that we can use \eqref{sth3.18} or \eqref{sth3.19} to generate the estimation of the unknown source $f$ when applying Algorithm 1. We show the numerical comparisons for using different estimation way in Algorithm 1 and Algorithm 2. And it can be seen that whatever approximation way, we almost get similar numerical result for all examples. In addition, it is obvious that the numerical effect is better for smaller noise level.  


\end{document}